\documentclass{amsart}
\usepackage{amsmath}
\usepackage{amssymb}
\usepackage{amsfonts}
\usepackage{amstext}
\theoremstyle{plain}\newtheorem{lem}[subsection]{Lemma}
\theoremstyle{plain}\newtheorem{df}[subsection]{Definition}
\theoremstyle{plain}\newtheorem{theo}[subsection]{\bf{Theorem}}
\theoremstyle{plain}\newtheorem{cor}[subsection]{\bf{Corollary}}
\newcommand{\bs}{\boldsymbol}
\newcommand{\R}{\mathbb{R}}

\newcommand{\C}{\mathbb{C}}
\newcommand{\Ha}{\mathbb{H}}
\def\ol{\overline}

\def\ls{\leqslant}
\def\Ra{\Rightarrow}
\def\ve{\varepsilon}

\def\fr{\forall}
\def\sm{\setminus}
\title[On boundary values of the
Cauchy-type integral] {On boundary values for rectifiable curves
of a generalization of the Cauchy-type integral related to the
Helmholtz operator in $\R^2$}
\author{Oleg F. Gerus}
\address{Oleg F. Gerus, Zhitomir State Pedagogical University,
Velyka Berdychivska str. 40, 10008 Zhitomir, Ukraine. {\it E-mail: ofg@com.zt.ua}}
\author{Michael Shapiro}
\address{Michael  Shapiro, Departamento de Matem\'aticas, ESFM del IPN, 07338 Mexico City,
 Mexico. {\it E-mail: shapiro@esfm.ipn.mx}}
\date{February 4, 2003}

\begin{document}
\begin{abstract}
There  are  considered  vector fields and  quaternionic
$\alpha$-hyperho\-lo\-mor\-phic functions in  a  domain of  $\R^2$
  which  generalize   the  notion  of  solenoidal  and  irrotational vector  fields.
There  are   established  sufficient conditions for the
corresponding Cauchy-type integral along a closed Jordan
rectifiable curve to be continuously extended onto the  closure of
a domain. The  Sokhotski-Plemelj-type formulas are proved as well.
\end{abstract}
\maketitle

\section{A generalization of holomorphy in
$\mathbb{R}^2$.}\label{s1}

\subsection{} Given a domain $\Omega\in\R^2$ consider a
$\C^3$-valued function $\bs{f}:\Omega\mapsto\C^3$. Let $\bs{i_1,
i_2, {\bs i_3}}$ be a canonical basis in $\mathbb{C}^3$, then any
$\bs{f}$ is of the form $\bs{f} = f_1\bs{i_1}+ f_2\bs{i_2}+
f_3\bs{i_3}$; moreover, we let $\mathbb{R}^2$ be a real linear
space with the  basis $\bs{i_1}, \bs{i_2}$. Of course, here $f_1,
f_2, f_3$ are complex-valued functions in $\Omega$ and
$\|\bs{f}\|^2:= |f_1|^2 + |f_2|^2 + |f_3|^2$.

\subsection{} In this paper we are interested in those vector-functions
$\bs{f}$ which are solutions of the following system:
\begin{equation}
\label{system1}
\begin{cases}
div\, \bs{f} = 0,\\ \bs{rot}\, \bs{f} = -\alpha \bs{f},
\end{cases}
\end{equation}

where $\alpha$ is a given complex number and for
$\bs{z}:=x\bs{i_1}+y\bs{i_2}\in\Omega$

\begin{eqnarray*}
          div \,\bs{f}&:=& \dfrac{\partial f_1}{\partial x} +
                                \dfrac{\partial f_2}{\partial y},\\
   \qquad  \bs{rot}\,\bs{f}&:=&\dfrac{\partial f_3}{\partial y}\bs{i_1} -
                                \dfrac{\partial f_3}{\partial x}\bs{i_2} +
                        \biggl( \dfrac{\partial f_2}{\partial x} -
                                \dfrac{\partial f_1}{\partial y}
                        \biggr) \bs{i_3}.
\end{eqnarray*}

Solutions of the system (\ref{system1}) are natural
generalizations of the notion of a solenoidal and irrotational
vector field which corresponds to $\alpha=0$. There are many
papers about both the case $\alpha=0$ and $\alpha \in
\mathbb{C}\setminus\{0\}$, see for instance \cite{gr}, \cite{krs},
\cite{roch}. Note  that  the  system  (\ref{system1}) can  be
naturally  considered  in  a  domain  of  $\R^3$  but  since  the
two--dimensional  case  has  its  essential  peculiarities,  we
present  here  the  latter  one,  and  the  former  will  be  done
elsewhere.
\subsection{} But it appears that there are deep mathematical reasons to
consider a more general system, namely, the following one:
\begin{equation}
\label{system2}
\begin{cases}
div\, \bs{f} = \alpha f_0,\\
\bs{rot}\, \bs{f} + \alpha \bs{f} = -
\bs{grad} \; f_0,
\end{cases}
\end{equation}
where $f_0$ is a $\mathbb{C}$-valued function and $
    \bs{grad}\,f_0:=\dfrac{\partial f_0}{\partial x}\bs{i_1}+
                         \dfrac{\partial f_0}{\partial y}\bs{i_2}.
$ Thus we shall be working with pairs ${\mathcal{F}}:= (f_0,
\bs{f})$ with norm $\|\mathcal{F}\|^2:=|f_0|^2+\|\bs{f}\|^2$ and
satisfying (\ref{system2}) from where certain conclusions will be
made for vector--functions satisfying (\ref{system1}).

\subsection{} Let $H^{(p)}_n$ be the Hankel functions of the kind $p\in\{1,2\}$
and of order $n\in\{ 0, 1, 2\}$ (see \cite{gr}), introduce the following notation:
\[
K_{\alpha,0}(\bs z):=
\begin{cases}
(-1)^p\dfrac{i\alpha}{4}H^{(p)}_0(\alpha\|\bs z\|),  &\text{if $\,\alpha\ne 0$},\\
0,                                               &\text{if
$\,\alpha=0$},
\end{cases}
\]
\[
\bs{K}_\alpha(\bs{z}):=
\begin{cases}
(-1)^p\dfrac{i\alpha}{4}H^{(p)}_1(\alpha\|\bs{z}\|)\dfrac{\bs{z}}{\|\bs{z}\|},
                                    &\text{if $\,\alpha\ne 0$},\\
-\dfrac{\bs{z}}{2\pi\|\bs{z}\|^2},    &\text{if $\,\alpha=0$},
\end{cases}
\]
where $\bs{z}\in\mathbb{R}^2\sm\{(0;0)\}$, $i$ is the complex imaginary unit in
$\mathbb{C}$, and
\begin{equation}\label{p}
p=\begin{cases}
      1, & \text{if $\,Im(\alpha)>0\,$ or $\,\alpha>0$},\\
      2, & \text{if $\,Im(\alpha)<0\,$ or $\,\alpha<0$}.\\
  \end{cases}
\end{equation}

Now the following pair:
\begin{equation}\label{kalv}
K_\alpha(\bs z):= \biggl(K_{\alpha,0}(\bs z),\, \bs{K}_\alpha(\bs z)\biggr)
\end{equation}
will play the role, in a sense, of the Cauchy kernel for the system (\ref{system2}).
It generates an analog of the Cauchy-type integral for a continuous pair
$\mathcal{F}=(f_0,\bs{f})$, $f_0:\Gamma\mapsto\mathbb{C}$,
$\bs{f}:\Gamma\mapsto\mathbb{C}^3$; by the formulas: if $\bs\zeta:= \xi\bs{i_1} +
\eta\bs{i_2}$, $\bs{\sigma}:= d\eta\bs{i_1} - d\xi\bs{i_2}$, and if $\Gamma$ is a
closed rectifiable Jordan curve which is the boundary of a bounded domain
$\Omega^+:=\Omega$, and $\Omega^-:= \R^2\setminus(\Omega^+\cup\Gamma)$ is its
complement;  then we define the pair
\begin{equation}\label{fiav}
  \Phi_\alpha[\mathcal{F}](\bs z):=
\biggl(\Phi_{\alpha,0}[\mathcal{F}](\bs z),
\bs{\Phi}_\alpha[\mathcal{F}](\bs z)\biggr),\qquad \bs z\in
\mathbb{R}^2\setminus\Gamma,
\end{equation}
with
\begin{equation}
\label{tres}
\begin{split}
\Phi_{\alpha,0}[\mathcal{F}](\bs z):=&-\int\limits_\Gamma
  \biggl(
      \bigl< \bs{K}_\alpha(\bs\zeta-\bs z),\bs{\sigma}
      \bigr> f_0(\bs\zeta) + \\
      & \qquad
       + \Bigl< \bigl[ \bs{K}_\alpha(\bs\zeta-\bs z),\bs{\sigma}
                 \bigr] +
               K_{\alpha,0}(\bs\zeta-\bs z) \cdot \bs{\sigma},\bs{f}(\bs\zeta)
          \Bigr>
  \biggr),
\end{split}
\end{equation}

\begin{equation}
\begin{split}
\bs\Phi_\alpha[\mathcal{F}](\bs z):=\int\limits_\Gamma &
   \biggl(
       \Bigl[
           \left[ \bs{K}_\alpha(\bs\zeta-\bs z),\bs{\sigma}
           \right]+ K_{\alpha,0}(\bs\zeta-\bs z)\cdot\bs\sigma, \bs{f}(\bs\zeta)
       \Bigr] -
    \biggr. \\
     \biggl.
    &
        -\Bigl<\bs K_\alpha(\bs\zeta-\bs z), \bs{\sigma}
         \Bigr>\bs{f}(\bs\zeta) +
    \biggr.\\
    \biggl.
    &
      + f_0(\bs\zeta)
      \Bigl(
         \bigl[ \bs{K}_\alpha(\bs\zeta-\bs z),\bs{\sigma}
         \bigr]
     + K_{\alpha,0}(\bs\zeta-\bs z)\cdot\bs{\sigma}
      \Bigr)
   \biggr),
\end{split}
\end{equation}
where $<\cdot,\cdot>$ and $[\cdot,\cdot]$ denote, respectively, the scalar and the
vector products in $\mathbb{C}^3$, i.e., for $\{\bs{a},\bs{b}\}$, $\bs{a}=
\sum_{k=1}^3 a_k \bs{i_k}$,
 $\bs{b}= \sum_{k=1}^3 b_k \bs{i_k}$,
$<\bs{a}, \bs{b}>:= \sum_{k=1}^3 a_k b_k$;
\[
[\bs{a}, \bs{b}]:= \begin{vmatrix}
                         \bs{i_1} & \bs{i_2} & \bs{i_3}\\
                         a_1   &    a_2   &   a_3   \\
                         b_1   &    b_2   &   b_3
                   \end{vmatrix}.
\]
We shall write $\Phi^+_\alpha[f]$ and $\Phi^-_\alpha[f]$ for the
respective restrictions of  $\Phi_\alpha[f]$ onto $\Omega^+$ and
$\Omega^-$.

\subsection{} If, in particular, the scalar component of the pair
$\mathcal{F}$ is identically zero, $\mathcal{F}=\{0; \bs f\}$,
this does not mean, in general, that $\Phi_\alpha[\mathcal{F}]$ is
vector-valued, which follows directly from (\ref{tres}).  Thus,
thinking of an analog of the Cauchy-type integral for the system
(\ref{system1}) as a purely vectorial object, we must exclude the
scalar component of
 $\Phi_\alpha[\mathcal{F}]$ in the following sense.  Let
$C(\Gamma;\mathbb{C}^3)$ be the complex linear space of all
$\mathbb{C}^3$-valued continuous vector--functions $\bs f$ on
$\Gamma$, introduce
\begin{equation}
 {\mathcal{M}}(\Gamma;\C^3):= C(\Gamma;\C^3) \cap
\Bigl\{
   \bs{f}:\int \limits_{\Gamma}
           \bigl<
              \left[
               \bs{K}_\alpha(\bs\zeta-\bs z), \bs{\sigma}
              \right]
             +
               K_{\alpha,0}(\bs\zeta-\bs z)\cdot\bs{\sigma},\bs{f}(\bs\zeta)
            \bigr> = 0,\ \bs z\not\in\Gamma
\Bigr\}.
\end{equation}
Now, for $\bs{f}\in {\mathcal{M}}(\Gamma; \mathbb{C}^3)$, the
analog of the Cauchy-type integral for the system (\ref{system1})
is given by the formula
\begin{equation*}
\begin{split}
\bs\Phi_\alpha[\bs f](\bs z):=\int\limits_\Gamma &
   \biggl(
       \Bigl[
           \left[ \bs{K}_\alpha(\bs\zeta-\bs z),\bs{\sigma}
           \right]+ K_{\alpha,0}(\bs\zeta-\bs z)\cdot\bs\sigma, \bs{f}(\bs\zeta)
       \Bigr] -
    \biggr. \\
     \biggl.
    &
        -\Bigl<\bs K_\alpha(\bs\zeta-\bs z), \bs{\sigma}
         \Bigr>\bs{f}(\bs\zeta)
   \biggr).
\end{split}
\end{equation*}
One may check up that now $\bs\Phi_\alpha[\bs{f}]$ is a solution
to (\ref{system1}).

\begin{theo}[Analogue of N. A. Davydov theorem (see \cite{dav}) for the system
(\ref{system2})] \label{theo6} Let $\Gamma$ be a closed
rectifiable Jordan curve, $f_0: \Gamma \mapsto \mathbb{C}$ and
$\bs{f}: \Gamma \mapsto \mathbb{C}^3$ be continuous functions,
${\mathcal{F}}:=(f_0,\bs{f})$, and let the integral
\begin{equation}
\Psi_\alpha[{\mathcal{F}}](\bs t):= \lim_{\delta\to 0}
\int\limits_{\Gamma\setminus \Gamma_{\bs t,\delta}}
\|K_\alpha(\bs\zeta-\bs t)\| \cdot \|\bs{\sigma}\| \cdot
\|{\mathcal{F}}(\bs\zeta) - {\mathcal{F}}(\bs t) \|, \qquad \bs
t\in\Gamma,
\end{equation}
where $\Gamma_{\bs t,\delta}:= \left\{\bs\zeta \in \Gamma :
\|\bs\zeta-\bs t\|\ls \delta\right\}$, exist uniformly with
respect to $\bs t\in \Gamma$. Then there exists the pair of
integrals
\[
      F_\alpha[\mathcal{F}](\bs t):= \biggl( F_{\alpha,0}[{\mathcal{F}}](\bs t),
                             \bs{F}_\alpha [{\mathcal{F}}](\bs t)\biggl),
      \qquad \bs t\in \Gamma,
\]
where
\begin{equation}\begin{split}
F_{\alpha,0}[{\mathcal{F}}](\bs t):=-\lim_{\delta \rightarrow
0}\int\limits_{\Gamma \setminus \Gamma _{\bs t,\delta }}&
 \bigl(\left<\bs{K}_\alpha(\bs\zeta-\bs t),\bs{\sigma}\right>(f_0(\bs\zeta)-f_0(\bs t))
\bigr.+\\
+&\left.\bigl< \left[\bs{K}_\alpha(\bs\zeta-\bs t),\bs{\sigma}\right]
+K_{\alpha,0}(\bs\zeta-\bs t)\cdot\bs{\sigma}, \left(\bs{ f}(\bs\zeta)-\bs{ f}(\bs
t)\right)\bigr>\right),
\end{split} \label{efv1} \end{equation}
\begin{equation}
\begin{split}
\bs{F}_\alpha[{\mathcal{F}}](\bs t):= \lim_{\delta \to 0}
\int\limits_{\Gamma\setminus\Gamma_{\bs t,\delta }} &\left(
 \bigl[
  \left[\bs{K}_\alpha(\bs\zeta-\bs t),\bs{\sigma}
  \right]+K_{\alpha,0}(\bs\zeta-\bs t)\cdot\bs{\sigma},
   \left(\bs{f}(\bs\zeta)-\bs{f}(\bs t)
   \right)
 \bigr]
\right.-\\
-& \left< \bs{K}_\alpha(\bs\zeta-\bs t),\bs{\sigma}
   \right>
    \left(\bs{f}(\bs\zeta)-\bs{f}(\bs t)
    \right)+\\
+&\bigl.(f_0(\bs\zeta)-f_0(\bs t))
    \left(
     \left[\bs{K}_\alpha(\bs\zeta-\bs t),\bs{\sigma}
     \right]+K_{\alpha,0}(\bs\zeta-\bs t)\cdot\bs{\sigma}
    \right)
  \bigr);
\end{split} \label{efv2} \end{equation}
moreover, the functions $\Phi_\alpha^\pm [{\mathcal{F}}]$ extend
continuously onto $\Gamma$ and the following analogues of the
Sokhotski-Plemelj formulas hold:
\begin{equation}
\begin{split}
\Phi_{\alpha,0}^{+}[{\mathcal{F}}](\bs
t)&=(I_{\alpha,\Gamma,0}(\bs t)+1)f_0(\bs t)-
\left<\bs{I}_{\alpha,\Gamma}(\bs t),\bs{f}(\bs t)\right>
+F_{\alpha,0}[{\mathcal{F}}](\bs t),\quad \bs t\in\Gamma,\\
\bs{\Phi}_{\alpha }^{+}[{\mathcal{F}}](\bs t)&=
\left[\bs{I}_{\alpha,\Gamma}(\bs t),\bs{f}(\bs t)\right] +
(I_{\alpha,\Gamma,0}(\bs t)+1)\bs{f}(\bs t)+ f_0(\bs
t)\bs{I}_{\alpha,\Gamma}(\bs t)
 +\bs{F}_{\alpha}[{\mathcal{F}}](\bs t), \ \bs t\in\Gamma,
\end{split}
\label{phiv+}
\end{equation}
\begin{equation}
\begin{split}
\Phi_{\alpha,0}^{-}[{\mathcal{F}}](\bs t)&=I_{\alpha,\Gamma,0}(\bs
t)f_0(\bs t) - \left<\bs{I}_{\alpha,\Gamma}(\bs t),\bs{f}(\bs
t)\right>
+F_{\alpha,0}[{\mathcal{F}}](\bs t),\quad \bs t\in\Gamma,\\
\bs{\Phi}_{\alpha }^{-}[{\mathcal{F}}](\bs t)&=
\left[\bs{I}_{\alpha,\Gamma}(\bs t),\bs{f}(\bs t)\right] +
I_{\alpha,\Gamma,0}(\bs t)\bs{f}(\bs t)+ f_0(\bs
t)\bs{I}_{\alpha,\Gamma}(\bs t)
 +\bs{F}_{\alpha}[\mathcal{F}](\bs t), \quad \bs t\in\Gamma,
\end{split}
\label{phiv-}
\end{equation}
where the following notation is used:\\
\begin{equation*}
I_{\alpha,\Gamma,0 }(\bs t):=-\alpha\iint\limits_{\Omega^{+}} K_{\alpha,0}(\bs\zeta
-\bs t) d\xi d\eta,
\end{equation*}
\begin{equation*}
\bs{I}_{\alpha ,\Gamma}(\bs t) :=-\alpha\iint\limits_{\Omega^{+}}
\bs{K}_\alpha(\bs\zeta-\bs t) d\xi d\eta,
\end{equation*}
and $\Phi^\pm_\alpha[{\mathcal{F}}](\bs t):=\lim\limits_{\Omega^\pm\ni \bs z\to \bs
t}\Phi_\alpha[{\mathcal{F}}](\bs z)$.
\end{theo}
The proof will be given after some preparatory work which is of an
independent interest.

\begin{cor}[Analogue of N. A. Davydov theorem for the system
(\ref{system1})]\label{17} Let $\Gamma$ be a closed rectifiable Jordan curve, $\bs
f\in{\mathcal{M}}(\Gamma;\C^3)$, and let the integral
\begin{equation}
\Psi_\alpha[\bs f](\bs t):= \lim_{\delta\to 0} \int\limits_{\Gamma\setminus
\Gamma_{\bs t,\delta}} \|K_\alpha(\bs\zeta-\bs t)\| \cdot \|\bs{\sigma}\| \cdot
\|\bs f(\bs\zeta) - \bs f(\bs t) \|
\end{equation}
exist uniformly with respect to $\bs t\in\Gamma$. Then there
exists the integral
\begin{equation}
\begin{split}
\bs{F}_\alpha \left[ \bs{f} \right](\bs t):= \lim_{\delta
\rightarrow0} \int\limits_{\Gamma\setminus\Gamma_{\bs t,\delta }}
\Bigl(&
 \bigl[
  \left[\bs{K}_\alpha(\bs\zeta-\bs t),\bs{\sigma}
  \right]+K_{\alpha,0}(\bs\zeta-\bs t)\cdot\bs\sigma,
   \left(\bs{f}(\bs\zeta)-\bs{f}(\bs t)
   \right)
 \bigr]
\Bigr.-
\\
&\Bigl. -\left< \bs{K}_\alpha(\bs\zeta-\bs t),\bs{\sigma}
 \right>
    \bigl(\bs{f}(\bs\zeta)-\bs{f}(\bs t)
    \bigr)
\Bigr); \label{efvc2}
\end{split}
\end{equation}
moreover, the functions $\bs{\Phi}_{\alpha }^{\pm}\left[\bs{f}\right](\bs z)$ extend
continuously onto $\Gamma$ and the following analogues of the Sokhotski-Plemelj
formulas hold:
\begin{equation}
\bs{\Phi}_{\alpha }^{+}\left[\bs{f}\right](\bs t)
=\left[\bs{I}_{\alpha,\Gamma}(\bs t),\bs{f}(\bs t)\right] +
(I_{\alpha,\Gamma,0}(\bs t)+1)\bs{f}(\bs t)
 +\bs{F}_{\alpha}\left[\bs{f}\right](\bs t), \quad \bs t\in\Gamma,
\label{phivc+}
\end{equation}
\begin{equation}
\bs{\Phi}_{\alpha }^{-}\left[\bs{f}\right](\bs t)=
\left[\bs{I}_{\alpha,\Gamma}(\bs t),\bs{f}(\bs t)\right] +
I_{\alpha,\Gamma,0}(\bs t)\bs{f}(\bs t)
 +\bs{F}_{\alpha}\left[\bs{f}\right](\bs t), \quad \bs t\in\Gamma,
\label{phivc-}
\end{equation}
where $\bs{\Phi}^\pm_\alpha\left[\bs{f}\right](\bs t):=
\lim\limits_{\Omega^\pm\ni \bs z\to \bs
t}\bs{\Phi}_\alpha\left[\bs{f}\right](\bs z)$.
\end{cor}

\section{Qua\-ter\-nions and qua\-ter\-nion-valued
$\alpha$-hyper\-holo\-morphic\\ functions in $\R^2$}

\subsection{} We shall denote, as usual, by $\Ha =
\Ha(\R)$ and $\Ha(\C)$ the sets of real and complex quaternions,
i.e., each quaternion is of the form $$a=\sum_{k=0}^3a_k\bs{i_k}$$
with $\{a_k\} \subset \R$ for real quaternions and $\{a_k\}
\subset \C$ for complex quaternions; ${\bs i_0} = 1 $ stands for
the unit and ${\bs i_1}, {\bs i_2}, {\bs i_3}$ stand for imaginary
units; the complex numbers imaginary unit in $\C$ will be denoted
by $i$. $\Ha$ has the structure of a real non-commutative,
associative algebra without zero divisors. $\Ha(\C)$ is a complex
non-commutative, associative algebra with zero divisors.

For a complex quaternion $a=\sum\limits_{k=0}^3a_k{\bs i_k}$ its
quaternionic conjugate is defined by $\ol
a:=a_0-\sum\limits_{k=1}^3a_k{\bs i_k}$. The module of a
quaternion $a$ coincides with its Euclidean norm : $|a| =
\|a\|_{\R^8} $. In particular, for $a \in \Ha$ we have $|a| =
\|a\|_{\R^4} $ and besides $|a|^2 = a \cdot \ol a = \ol a \cdot a$
while for a complex quaternion $|a|^2 \neq a \cdot \ol a $. What
is more, for $a,b$ from $\Ha$ there holds: $|a \cdot b| = |a|
\cdot |b|$ which is extremely important working with real
quaternions. For complex quaternions the situation is different.
\begin{lem}
\label{l1} $|ab|\ls \sqrt{2}\cdot |a|\cdot |b|$ for all $a,b\in %
\Ha(\C)$.
\end{lem}
\begin{proof}
Let $$a:=\sum_{k=0}^3a_k{\bs i_k},\
a_k:=\alpha_k+i\lambda_k,\qquad \alpha_k,\lambda_k\in\R,$$
$$b:=\sum_{k=0}^3b_k{\bs i_k},\ b_k:=\beta_k+i\gamma_k,\qquad
\beta_k,\gamma_k\in\R.$$

We have $a=a'+ia''$ and $b=b'+ib''$, where $a',a'',b',b''$ are
real quaternions. Since
$$|a|^2=\sum\limits_{k=0}^{3}(\alpha_k^2+\lambda_k^2)= |a'|^2+
|a''|^2,\quad |b|^2=\sum\limits_{k=0}^{3}(\beta_k^2+\gamma_k^2)=
|b'|^2+ |b''|^2,$$ then $$ (|a|\cdot|b|)^2=|a'|^2|b'|^2+
|a'|^2|b''|^2+ |a''|^2|b'|^2+ |a''|^2|b''|^2. $$ Therefore
\begin{equation}
\begin{split}
|ab|^2&=|a'b'-a''b''+i(a'b''+a''b')|^2=
|a'b'-a''b''|^2+|a'b''+a''b'|^2=\\
&=\ol{(a'b'-a''b'')}\cdot(a'b'-a''b'')+
\ol{(a'b''+a''b')}\cdot(a'b''+a''b')=\\
&=(\ol{a'b'}-\ol{a''b''})\cdot(a'b'-a''b'')+
(\ol{a'b''}+\ol{a''b'})\cdot(a'b''+a''b') =(|a|\cdot|b|)^2 + d,
\end{split}
\label{mab}
\end{equation}
where $\Ha\ni d=\ol{a'b''}a''b'+\ol{a''b'}a'b''-\ol{a'b'}a''b''-\ol{a''b''}a'b',$
and
\begin{equation}
|d|\ls2|a'b'|\cdot|a''b''|+2|a'b''|\cdot|a''b'|\ls(|a|\cdot|b|)^2.
\label{d}
\end{equation}
Combining (\ref{mab}) and (\ref{d}), we obtain the assertion of
Lemma.
\end{proof}

\subsection{} Let $\Omega$ be a domain in the
plane $\R^2$, we consider $\Ha(\C)$-valued functions defined in
the domain $\Omega$. On the bi-$\Ha(\C)$-module
$C^2(\Omega;\Ha(\C))$ there introduced the two-dimensional
Helmholtz operator with a wave number $\lambda \in \C$:
$$_\lambda\Delta:=\Delta_{\R^2}+\hspace{0pt}^\lambda M,$$
where $\Delta_{\R^2}=\partial_{1}^2 + \partial_{2}^2$, $\partial_k =
\frac{\partial}{\partial x_k} $ and for $a \in \Ha(\C)$ we denote $^\lambda M$ the
operator of the multiplication by $\lambda$ on the left-hand side, analogously for
$M^\lambda$. The operators $\ol \partial :=
\partial_1 + i\partial_2 $ and $\partial :=
\partial_1 - i\partial_2$ determine, respectively, classes of
holomorphic and anti-holomorphic functions of the complex
variable, and the following factorization holds: $$\partial \circ
\ol \partial = \ol \partial \circ \partial = \Delta_{\R^2}.$$
Consider the following partial differential operators with
quaternionic coefficients: $$_{st}\partial := {\bs i_1} \cdot
\partial_1 + {\bs i_2} \cdot
\partial_2;\
_{st} \ol \partial := \ol{{\bs i_1}} \cdot \partial_1 + \ol{{\bs
i_2}} \cdot
\partial_2;$$
$$\partial_{st} := \partial_1 \circ M^{{\bs i_1}} + \partial_2
\circ M^{{\bs i_2}};\ \ol \partial_{st} := \partial_1 \circ M^{\ol
{\bs i_1}} +
\partial_2 \circ M^{\ol {\bs i_2}}.$$
 The following equalities can be easily verified:
$$\partial_{st}
\circ \ol \partial_{st} = \ol \partial_{st} \circ
\partial_{st} = \Delta_{\R^2} = \hspace{0pt}_{st}\partial \circ \hspace{0pt}_{st} \ol \partial =
 \hspace{0pt}_{st} \ol \partial \circ \hspace{0pt}_{st}\partial,$$ which mean that
$$_{st}\partial^2 =
\partial^2_{st}=- \Delta_{\R^2}.$$
For $\alpha\in\C$ to be a complex square root of $\lambda\in\C$, i.e. $\alpha^2 =
\lambda$, set
$$_\alpha\partial :=
\partial_{st}+ ^{\alpha}\!\!M;\ \partial_{\alpha}:=\hspace{0pt}_{st}\partial +
M^\alpha.$$ Then we have the following factorizations of the Helmholtz operator:
$$_\lambda\Delta=-\partial_\alpha \circ
\partial_{-\alpha}=-\partial_{-\alpha} \circ \partial_\alpha=-_\alpha\partial \circ
_{-\alpha}\partial=-_{-\alpha}\partial \circ _\alpha\!\partial.$$
In analogy with the usual notion of a holomorphic function,
consider the following definition of $\alpha$-hyperholomorphic
functions in $\R^2$.

\begin{df}[\cite{MS97}]\label{df} Let $f \in
C^1(\Omega,\Ha(\C))$, $f$ is called $\alpha$-hyperholomorphic if $_{\alpha}\partial
f \equiv 0$ in $\Omega$.
\end{df}

Of course, more exactly such functions should be called, for instance,
left-$\alpha$-hyperholomorphic because there is a "symmetric" definition for
$\partial_{\alpha}$, as well as for $_{\alpha}\ol
\partial$ and $\ol \partial_{\alpha}$. We shall deal with
the above case only.

Such a definition for $\alpha$-hyperholomorphic functions was introduced in
\cite{MS97} both for complex and quaternionic values of $\alpha$, and some essential
properties were established there. Main integral formulas for
$\alpha$-hyperholomorphic functions were constructed in \cite{MS98}. All proofs and
details can be found in these papers, see also \cite[Appendix 4]{krs}. Some
developments of the topic are presented in \cite{roch} and \cite{GSS}. One can find
much more relevant bibliographical references in all these papers.

\subsection{} In what follows we shall be in
need of some properties of the Hankel functions $H_n^{(p)}(t)$,
$t\in\C$, (see \cite{gr}) which we concentrate in this section for
the reader's convenience. The following equalities are valid:
\begin{equation}\label{dh1}
\frac{d}{dt}H^{(p)}_1(t)=\frac{1}{2}\left(H^{(p)}_0(t)-
H^{(p)}_2(t)\right),
\end{equation}
\begin{equation}\label{dh0}
\frac{d}{dt}H^{(p)}_0(t)=-H^{(p)}_1(t),
\end{equation}
\begin{equation}\label{tha2}
tH^{(p)}_2(t)=2H^{(p)}_1(t)-tH^{(p)}_0(t),
\end{equation}
 and the following series expansions of the
Hankel functions $H_0^{(p)}(t)$ and $H_1^{(p)}(t)$ hold:
\begin{equation}\label{ha0}
H_0^{(p)}(t)=
\left(1-(-1)^p\frac{2i}{\pi}(\log\frac{t}{2}+{\boldsymbol
C})\right) \sum\limits_{k=0}^{\infty} \frac{(-1)^k
t^{2k}}{2^{2k}(k!)^2} + \frac{2i}{\pi} \sum\limits_{k=1}^{\infty}
\frac{(-1)^{k+p} t^{2k}}{2^{2k}(k!)^2}\sum\limits_{m=1}^{k}\frac{1}{m}, \\
\end{equation}
\begin{equation}\label{ha1}
\begin{split}
    H_1^{(p)}(t)= & \left(1-(-1)^p\frac{2i}{\pi}(\log\frac{t}{2}+{\boldsymbol
C})\right) \sum\limits_{k=0}^{\infty} \frac{(-1)^k
t^{2k+1}}{2^{2k+1}k!(k+1)!} +(-1)^p \left(\frac{2i}{\pi t} +
\frac{i t}{2\pi}\right)+\\
        &+ \frac{i}{\pi} \sum\limits_{k=1}^{\infty} \frac{(-1)^{k+p}
t^{2k+1}}{2^{2k+1}k!(k+1)!}\left(\sum\limits_{m=1}^{k+1}
\frac{1}{m} + \sum\limits_{m=1}^{k} \frac{1}{m}\right),
\end{split}
\end{equation}
where $\boldsymbol C$ is the Euler constant.

\section{Quaternionic generalization of the Cauchy-type
integral}\label{sec}

\subsection{} Given $\alpha \in \C$ and real
quaternions $z:=x {\bs i_1}+y {\bs i_2}$, $\zeta:=\xi {\bs
i_1}+\eta {\bs i_2}$ thinkable as points of the Euclidean space
$\R^2$ equipped with the additional structure of quaternionic
multiplication; introduce the notation: $$\theta_{\alpha}(z):=
\begin{cases}
(-1)^p\dfrac{i}{4}H_0^{(p)}(\alpha|z|),& \text{if\, $\alpha\not=0$},\\
\dfrac1{2\pi}\log |z|,& \text{if\, $\alpha=0$},
\end{cases}$$
where $p$ depends on $\alpha$ by the formula (\ref{p}). It is a well known fact (see
e.g. \cite{v71}) that the function $\theta_{\alpha}$ is the fundamental solution of
the Helmholtz operator
$\Delta_{\alpha^2}:=\Delta_{\R^2}+\hspace{0pt}^{\alpha^2}\!M$, written down for all
values of $\alpha$.

The $\alpha$-hyperholomorphic Cauchy kernel, i.e., the fundamental solution to
o\-pe\-rators $_{\alpha}\partial$ and $\partial_{\alpha}$, is defined as
$$K_\alpha(z):=-_{-\alpha}\partial[\theta_\alpha](z)=
-\partial_{-\alpha}[\theta_\alpha](z).$$

Hence one has explicitly:
\begin{equation}\label{kal}
    K_\alpha(z)=
\begin{cases}
(-1)^p\dfrac{i\alpha}{4}\left(H_1^{(p)}(\alpha|z|) \dfrac{z}{|z|}
 + H_0^{(p)}(\alpha |z|) \right), &\text{if\, $\alpha\not=0$},\\
-\dfrac{z}{2\pi|z|^2},
 &\text{if\, $\alpha=0$}.
 \end{cases}\end{equation}
Now, for a continuous function $f:\Gamma\to\Ha(\C)$ and
$\sigma:=d\eta {\bs i_1}-d\xi {\bs i_2}$ the Cauchy-type integral
of $f$ is given by the formula
\begin{equation}\label{fia}
    \Phi_\alpha[f](z):=\int_\Gamma
K_\alpha(\zeta-z)\cdot\sigma\cdot f(\zeta),\qquad z\in
\R^2\sm\Gamma;
\end{equation} as in the previous section here $\Gamma$ is a
closed rectifiable Jordan curve.

\begin{theo} \label{mait}
Let $\Gamma $ be a closed rectifiable Jordan curve, $f:\Gamma \rightarrow %
\mathbb{H}(\C)$ be a continuous function, and let the integral
\begin{equation}
\Psi _{\alpha }\lbrack f\rbrack (t):=\lim_{\delta \rightarrow
0}\int\limits_{\Gamma \setminus \Gamma _{t,\delta }}|K_{\alpha }(\zeta
-t)|\cdot |\sigma |\cdot |f(\zeta )-f(t)|,\qquad t\in \Gamma ,  \label{psi}
\end{equation}
where $\Gamma _{t,\delta }:=\{\zeta \in \Gamma :|\zeta -t|\ls
\delta \} $, exist uniformly with respect to $t\in \Gamma $. Then
there exists the integral
\begin{equation}
F_{\alpha }\lbrack f\rbrack (t):=\lim_{\delta \rightarrow
0}\int\limits_{\Gamma \setminus \Gamma _{t,\delta }}K_{\alpha
}(\zeta -t)\cdot \sigma \cdot (f(\zeta )-f(t)),\qquad t\in \Gamma;
\label{ef}
\end{equation}
moreover, the functions $\Phi _{\alpha }^{\pm}\lbrack f\rbrack $ extend continuously
onto $\Gamma$, and the following analogues of the Sokhotski-Plemelj formulas hold:
\begin{equation}
\Phi _{\alpha }^{+}\lbrack f\rbrack (t)=(I_{\alpha ,\Gamma
}(t)+1)f(t)+F_{\alpha }\lbrack f\rbrack (t),\quad t\in \Gamma ,  \label{phi+}
\end{equation}
\begin{equation}
\Phi _{\alpha }^{-}\lbrack f\rbrack (t)=I_{\alpha ,\Gamma }(t)f(t)+F_{\alpha
}\lbrack f\rbrack (t),\quad t\in \Gamma ,  \label{phi-}
\end{equation}
where $\Phi^\pm_\alpha[f](t):=\lim\limits_{\Omega^\pm\ni z\to t}\Phi_\alpha[f](z)$,
and
\begin{equation*}
    I_{\alpha ,\Gamma }(t):=-\alpha \iint\limits_{\Omega^{+}} K_\alpha(\zeta-t)d\xi d\eta.
\end{equation*}
\end{theo}

The proof is based on several lemmas which are of interest by
themselves.

\begin{lem}
\label{l2} The limit (\ref{ef}) exists uniformly with respect to $t\in
\Gamma $ and $F_{\alpha }\lbrack f\rbrack $ is a continuous function on $%
\Gamma $.
\end{lem}

\begin{lem}
\label{l3}
\begin{equation}
\Phi _{\alpha }\lbrack 1\rbrack (z)=\begin{cases}
I_{\alpha,\Gamma}(z)+1, &z\in \Omega^+,\cr I_{\alpha,\Gamma}(z),
&z\in \Omega^-. \end{cases}  \label{fi}
\end{equation}
\end{lem}

\begin{lem}
\label{l4} $I_{\alpha ,\Gamma }$ is a continuous function in
$\mathbb{R}^{2}$.
\end{lem}

\section{Proof of the results of Section \ref{sec}. }

\stepcounter{subsection}
\begin{proof}[{\rm\thesubsection.} Proof of Lemma \ref{l2}] Denote
$$\Psi_\alpha(\delta,t):=
\int\limits_{\Gamma\sm\Gamma_{t,\delta}}
|K_\alpha(\zeta-t)|\cdot|\sigma|\cdot|f(\zeta)-f(t)|.$$

$$F_\alpha(\delta,t):=
\int\limits_{\Gamma\sm\Gamma_{t,\delta}}
K_\alpha(\zeta-t)\cdot\sigma\cdot(f(\zeta)-f(t)).$$

We have
$F_\alpha(\delta,t)=\sum\limits_{k=0}^{3}(F^{(1)}_{\alpha,k}(\delta,t)
+ iF^{(2)}_{\alpha,k}(\delta,t)){\bs i_k}$, where
$F^{(1)}_{\alpha,k}$ and $F^{(2)}_{\alpha,k}$ are real-valued
functions.

Under conditions of Theorem \ref{mait} the function
$\Psi_\alpha(\delta,t)$ tends to the finite limit
$\Psi_\alpha(t)$, when $\delta\to0$, uniformly with respect to
$t\in\Gamma$.  Using the criterion of uniform convergence for the
integral and Lemma \ref{l1}, we get that  for   $\fr\ve>0\
\exists\delta(\ve)>0 \ \fr t\in\Gamma:$
\begin{equation*}
\begin{split}
0<\delta_1&<\delta_2 <\delta(\ve)\Ra \\
&\Ra\Psi_\alpha(\delta_1,t)-\Psi_\alpha(\delta_2,t)=
\int\limits_{\Gamma_{t,\delta_2}\sm\Gamma_{t,\delta_1}}
|K_\alpha(\zeta-t)|\cdot|\sigma|\cdot|f(\zeta)-f(t)| < \ve\Ra\\
&\Ra\left|F_\alpha(\delta_1,t)-F_\alpha(\delta_2,t)\right| =
\left|\ \int\limits_{\Gamma_{t,\delta_2}\sm\Gamma_{t,\delta_1}}
K_\alpha(\zeta-t)\cdot\sigma\cdot(f(\zeta)-f(t)) \right|\ls\\
&\hspace{15pt} \ls 2 \int \limits_{\Gamma_{t,\delta_2} \sm
\Gamma_{t,\delta_1}}
|K_\alpha(\zeta-t)|\cdot|\sigma|\cdot|f(\zeta)-f(t)|
 < 2\ve\Ra\\
\end{split}
\end{equation*}
\begin{equation}
\Ra
\left|F^{(j)}_{\alpha,k}(\delta_1,t)-F^{(j)}_{\alpha,k}(\delta_2,t)\right|<2\ve\quad
(j=1,2;\ k=0,\dots,3). \label{phi}
\end{equation}

Therefore for each
fixed $t\in\Gamma$ there exist limits
$$F^{(j)}_{\alpha,k}(t):=\lim\limits_{\delta\to0}F^{(j)}_{\alpha,k}(\delta,t)\quad
(j=1,2;\ k=0,\dots,3)
$$
and consequently there exists
\begin{equation}
F_\alpha[f](t)=\lim\limits_{\delta\to0}F_\alpha(\delta,t).
\label{lim} \end{equation}

Proceeding to the limit as $\delta_1\to0$ in inequality
(\ref{phi}) we obtain  that  for  $\fr\ve>0\ \exists\delta(\ve)>0\
\fr t\in\Gamma$:
\begin{equation*}
\begin{split} 0<\delta<\delta(\ve) &\Ra
\left|F^{(j)}_{\alpha,k}(\delta,t)-F^{(j)}_{\alpha,k}(t)\right| \ls 2\ve\quad
(j=1,2;\
k=0,\dots,3)\Ra\\
&\Ra\left|F_\alpha(\delta,t)-F_\alpha[f](t)\right| \ls 4\sqrt2\ve,
\end{split}
\label{phi1}
\end{equation*}
and this is all.
\end{proof}

\stepcounter{subsection}
\begin{proof}[{\rm\thesubsection.} Proof of Lemma \ref{l3}]
Let $\alpha\neq0$. We have
\begin{equation}
\Phi_\alpha[1](z)=\int\limits_{\Gamma}K_\alpha(\zeta-z)\cdot\sigma
= (-1)^{p-1}\frac{i}{4}\alpha({\bs i_1}+I_2{\bs i_3}-I_3),
\label{27}
\end{equation}
where
\begin{align*}
I_1=&\int\limits_{\Gamma}\frac{H_1^{(p)}(\alpha|\zeta-z|)}{|\zeta-z|}
\left((\xi-x)d\eta-(\eta-y)d\xi \right),\\
I_2=&\int\limits_{\Gamma}\frac{H_1^{(p)}(\alpha|\zeta-z|)}{|\zeta-z|}
\left((\xi-x)d\xi+(\eta-y)d\eta \right), \\
I_3=&\int\limits_{\Gamma}H_0^{(p)}(\alpha|\zeta-z|) (d\eta {\bs
i_1} -d\xi {\bs i_2}).
\end{align*}

Let $z\in \Omega^+$, let $\rho>0$ be such that
$B(z,\rho):=\{\zeta\in\C: |\zeta-z|\ls\rho\}$ is contained in
$\Omega^+$, and let $\gamma_\rho$ be the boundary of $B(z,\rho)$.

Using the Green formula and the equalities (\ref{dh1}), (\ref{tha2}) we get:
\begin{align*}
  &\int\limits_{\Gamma-\gamma_\rho}  \frac{H_1^{(p)}(\alpha|\zeta-z|)}
  {|\zeta-z|} \left((\xi-x)d\eta-(\eta-y)d\xi \right)= \\ &
=\iint\limits_{\Omega^+\sm B(z,\rho)}\!\!\left(\frac{\partial }{\partial\xi }\left(
\frac{H_1^{(p)}(\alpha|\zeta-z|)} {|\zeta-z|} (\xi-x)\right) + \frac{\partial
}{\partial\eta}\left(\frac{H_1^{(p)}(\alpha|\zeta-z|)} {|\zeta-z|}
(\eta-y)\right)\right)d\xi\,d\eta =\\
& =\frac{1}{2}\iint\limits_{\Omega^+\sm B(z,\rho)} \left(\alpha
H_0^{(p)}(\alpha|\zeta-z|) - \alpha H_2^{(p)}(\alpha|\zeta-z|) +
2\frac{H_1^{(p)}(\alpha|\zeta-z|)} {|\zeta-z|}\right)d\xi\,d\eta=\\
&=\iint\limits_{\Omega^+\sm B(z,\rho)}\alpha H_0^{(p)}(\alpha|\zeta-z|)d\xi\,d\eta,
\end{align*}
\begin{align*}
  \int\limits_{\gamma_\rho}  \frac{H_1^{(p)}(\alpha|\zeta-z|)}
  {|\zeta-z|} \left((\xi-x)d\eta-(\eta-y)d\xi \right)  =& 2\pi\rho
H_1^{(p)}(\alpha\rho)=\\ =& (-1)^p\frac{4i}{\alpha}+o(1)\ \mbox{as}\ \rho\to0.
\end{align*}
Hence
\begin{equation} \begin{split} I_1 &
=\lim_{\rho\to0}\left(\int\limits_{\Gamma-\gamma_\rho} +
\int\limits_{\gamma_\rho} \right)
\frac{H_1^{(p)}(\alpha|\zeta-z|)}
  {|\zeta-z|} ((\xi-x)d\eta-(\eta-y)d\xi)=\\ &
=\iint\limits_{\Omega^+} \alpha H_0^{(p)}(\alpha|\zeta-z|)d\xi\,d\eta
+(-1)^p\frac{4i}{\alpha}.
\end{split}
\label{28a}
\end{equation}
Furthermore,
\begin{align*}
  &\int\limits_{\Gamma-\gamma_\rho}
  \frac{H_1^{(p)}(\alpha|\zeta-z|)}{|\zeta-z|}
\left((\xi-x)d\xi+(\eta-y)d\eta \right) =  \\
& =\iint\limits_{\Omega^+\sm B(z,\rho)}\!\! \left(\frac{\partial }{\partial\xi
}\left( \frac{H_1^{(p)}(\alpha|\zeta-z|)} {|\zeta-z|}(\eta-y) \right) -
\frac{\partial }{\partial\eta}\left(\frac{H_1^{(p)}(\alpha|\zeta-z|)}
{|\zeta-z|}(\xi-x)
\right)\right)d\xi\,d\eta =\\
& =\iint\limits_{\Omega^+\sm B(z,\rho)} \frac{\partial }{\partial
|\zeta-z|}\left(\frac{H_1^{(p)}(\alpha|\zeta-z|)} {|\zeta-z|} \right)\!\!
\left(\frac{\partial |\zeta-z|}{\partial \xi}(\eta-y) - \frac{\partial |\zeta-z|}
{\partial \eta}(\xi-x)\right)d\xi\,d\eta=\\&=0,
\end{align*}
\begin{align*}
  \int\limits_{\gamma_\rho}
  \frac{H_1^{(p)}(\alpha|\zeta-z|)}{|\zeta-z|}
\left((\xi-x)d\xi+(\eta-y)d\eta \right)=0
\end{align*}
and, consequently,
\begin{equation}
I_2=0.
\label{28}
\end{equation}

Analogously, using the equality (\ref{dh0}), we have
\begin{align*}
\int\limits_{\Gamma-\gamma_\rho} & H_0^{(p)}(\alpha|\zeta-z|)
(d\eta {\bs i_1}-d\xi
 {\bs i_2})=\\
& =\iint\limits_{\Omega^+\sm B(z,\rho)}\left(\frac{\partial
H_0^{(p)} (\alpha|\zeta-z|)} {\partial\xi} {\bs i_1} +
\frac{\partial H_0^{(p)} (\alpha|\zeta-z|)} {\partial\eta} {\bs
i_2}\right)d\xi\,d\eta=\\ & =-\iint\limits_{\Omega^+\sm B(z,\rho)}
\frac{\alpha H_1^{(p)}(\alpha|\zeta-z|)}{|\zeta-z|}
(\zeta-z)d\xi\,d\eta,
\end{align*}
\begin{align*}
\int\limits_{\gamma_\rho}H_0^{(p)}(\alpha|\zeta-z|) (d\eta {\bs
i_1}-d\xi {\bs i_2})=0
\end{align*}
and
\begin{equation}
I_3=-\iint\limits_{\Omega^+} \frac{\alpha
H_1^{(p)}(\alpha|\zeta-z|)}{|\zeta-z|} (\zeta-z)d\xi\,d\eta.
\label{29}
\end{equation}

Thus, from (\ref{27}) -- (\ref{29}) we have
\begin{align*}
\Phi_\alpha[1](z)= 1&+(-1)^{p-1}\frac{i\alpha^2}{4} \iint\limits_{\Omega^+} \left(
H_0^{(p)}(\alpha|\zeta-z|) \right.+ \\
& + \left.\frac{H_1^{(p)}(\alpha|\zeta-z|)}{|\zeta-z|}
(\zeta-z)\right)d\xi\,d\eta=I_{\alpha,\Gamma}(z)+1.
\end{align*}

Let now $\alpha=0$. We have
\begin{equation*}
\begin{split}
    \Phi_0[1](z) &= \int\limits_{\Gamma}K_0(\zeta-z) \cdot\sigma =
     -\frac1{2\pi}\int\limits_{\Gamma}\frac{\zeta-z}{|\zeta-z|^2}\cdot\sigma=\\
        &=\frac1{2\pi}\int\limits_{\Gamma}\frac{(\xi-x)d\eta-(\eta-y)d\xi}{|\zeta-z|^2}
        +\frac1{2\pi}\int\limits_{\Gamma} \frac{(\xi-x)d\xi+
        (\eta-y)d\eta}{|\zeta-z|^2}{\bs i_3}.
\end{split}
\end{equation*}
Going on along the same way as in the computation of the integrals
$I_1$ and $I_2$, and using the Green formula, we obtain that
$\Phi_0[1](z)=1$.

In the case of $z\in \Omega^-$ the proof of (\ref{fi}) is
simplified because of the continuity of the kernel $K_\alpha$ on
$\Omega^+$.
\end{proof}

\stepcounter{subsection}
\begin{proof}[{\rm\thesubsection.} Proof of Lemma \ref{l4}]
Making use of the series expansions of the Hankel functions (\ref{ha0}), (\ref{ha1})
we obtain:
\begin{equation}
I_{\alpha,\Gamma}(z)=\frac{i\alpha}{8}
\left((-1)^{p-1}I_{\alpha,\Gamma}^{(1)}(z) +
\frac{2i}{\pi}I_{\alpha,\Gamma}^{(2)}(z) + (-1)^p\frac{4i}{\pi}
I_\Gamma^{(3)}(z)\right), \label{l41}
\end{equation}
where
\begin{align*}
I_{\alpha,\Gamma}^{(1)}(z) :=&\iint\limits_{\Omega^+}
\sum\limits_{k=0}^{\infty}\alpha^{2k+1}|\zeta-z|^{2k} (a_{k,p} +
b_{k,p}\alpha(\zeta-z))d\xi\,d\eta,\\
I_{\alpha,\Gamma}^{(2)}(z) :=&\iint\limits_{\Omega^+}
\log|\zeta-z| \sum\limits_{k=0}^{\infty}(-1)^k
\frac{\alpha^{2k+1}|\zeta-z|^{2k}}{2^{2k}k!(k+1)!} (2(k+1) +
\alpha(\zeta-z))d\xi\,d\eta,\\
I_\Gamma^{(3)}(z) :=&\iint\limits_{\Omega^+}
\frac{\zeta-z}{|\zeta-z|^2}d\xi\,d\eta,
\end{align*}
and $a_{k,p}$, $b_{k,p}$ are complex coefficients.

The continuity of $I_{\alpha,\Gamma}^{(1)}$ follows from the
continuity of the integrand.

Let us prove the continuity of $I_\Gamma^{(3)}$. For an arbitrary $z\in\C$ and a
measurable $E\subset\C$ set
$$I_E(z):=\iint\limits_{E}
\frac{\zeta-z}{|\zeta-z|^2}d\xi\,d\eta.$$

Let us fix any point $z_0\in\C$. For an arbitrary $z\in\C$ we have
\begin{align*}
I_\Gamma^{(3)}(z_0)-I_\Gamma^{(3)}(z)&=I_{\Omega^+\cap B(z_0,\rho)}(z_0)+
I_{(\Omega^+\sm B(z_0,\rho))\cap B(z,\rho)}(z_0) +\\&+ I_{\Omega^+\sm (B(z,\rho)\cup
B(z_0,\rho))}(z_0) - I_{\Omega^+\cap B(z_0,\rho)}(z) -\\&- I_{(\Omega^+\sm
B(z_0,\rho))\cap B(z,\rho)}(z) - I_{\Omega^+\sm (B(z,\rho)\cup B(z_0,\rho))}(z).
\end{align*}

Fix an arbitrary $\ve>0$. Since $|I_{E\cap
B(z_1,\rho)}(z_2)|\ls16\rho$ for an arbitrary $\rho>0$,
$z_1\in\C$, $z_2\in\C$, $E\subset\C$, there exists $\rho(\ve)>0$
such that $|I_{E\cap B(z_1,\rho)}(z_2)|\ls\frac{\ve}{6}$.
Therefore
\begin{align*}
\left|I_\Gamma^{(3)}(z_0)-I_\Gamma^{(3)}(z)\right| &
\ls\frac{2\ve}{3}+ |I_{\Omega^+\sm (B(z,\rho)\cup
B(z_0,\rho))}(z_0)- I_{\Omega^+\sm (B(z,\rho)\cup B(z_0,\rho))}(z)|\ls\\
&\ls\frac{2\ve}{3} + \frac{4}{\pi}|z_0-z|\iint
\limits_{\Omega^+\sm (B(z,\rho)\cup B(z_0,\rho))}
\frac{d\xi\,d\eta}{|\zeta-z_0||\zeta-z|}.
\end{align*}
Under the condition $|z_0-z|<\frac{\rho(\ve)}{2}$ we get:
$$\iint \limits_{\Omega^+\sm (B(z,\rho)\cup
B(z_0,\rho))} \frac{d\xi\,d\eta}{|\zeta-z_0||\zeta-z|}\ls
4\pi\log\frac{d}{\rho(\ve)},$$ where
$d=\max\limits_{t\in\Gamma}|z_0-t|$. By choosing
$|z_0-z|<\min\left\{\frac{\rho(\ve)}{2};
\ve\left(48\log\frac{d}{\rho(\ve)}\right)^{-1}\right\}$ we obtain
$$\left|I_\Gamma^{(3)}(z_0)-I_\Gamma^{(3)}(z)\right|<\ve.$$

To prove the continuity of $I^{(2)}_{\alpha,\Gamma}(z)$ fix any
point $z_0\in \C$. For any $z\in\C$ denote $\delta:=|z-z_0|$,
$\Omega_1^+:=B(z_0,3\delta)\cap \Omega^+$,
$\Omega_2^+:=\Omega^+\sm \Omega_1^+$.

Denote
$$R(\zeta,z):=\log|\zeta-z|
\sum\limits_{k=0}^{\infty}(-1)^k
\frac{\alpha^{2k+1}|\zeta-z|^{2k}}{2^{2k}k!(k+1)!} (2(k+1) +
\alpha(\zeta-z)).$$ Then
\begin{align*}
|I_{\alpha,\Gamma}^{(2)}(z)-I_{\alpha,\Gamma}^{(2)}(z_0)|\ls&
\left|\iint\limits_{\Omega_1^+}R(\zeta,z_0)d\xi\,d\eta\right| +
\left|\iint\limits_{\Omega_1^+}R(\zeta,z)d\xi\,d\eta\right|+\\
+& \left|\iint\limits_{\Omega_2^+}\left(R(\zeta,z)-R(\zeta,z_0)\right)d\xi\,d\eta
\right|=: I_4+I_5+I_6,
\end{align*}
and setting $4\delta<1$ we get
\begin{align*}
I_4\ls&\iint\limits_{\Omega_1^+}|\log|\zeta-z_0|| \sum\limits_{k=0}^{\infty}
\frac{|\alpha|^{2k+1}|\zeta-z_0|^{2k}}{2^{2k}k!(k+1)!} (2(k+1) +
|\alpha|\cdot|\zeta-z_0|)d\xi\,d\eta\ls\\
\ls& \sum\limits_{k=0}^{\infty}
\frac{|\alpha|^{2k+1}|3\delta|^{2k}}{2^{2k}k!(k+1)!} (2(k+1) +
|\alpha|\cdot3\delta)2\pi \int\limits_{0}^{3\delta}\rho\log\frac{1}{\rho}\,d\rho =
o(1)\ as\ \delta\to0,
\end{align*}
\begin{align*}
I_5\ls&\iint\limits_{\Omega_1^+}|\log|\zeta-z|| \sum\limits_{k=0}^{\infty}
\frac{|\alpha|^{2k+1}|\zeta-z|^{2k}}{2^{2k}k!(k+1)!} (2(k+1) +
|\alpha|\cdot|\zeta-z|)d\xi\,d\eta\ls\\
\ls& \sum\limits_{k=0}^{\infty}
\frac{|\alpha|^{2k+1}|4\delta|^{2k}}{2^{2k}k!(k+1)!} (2(k+1) +
|\alpha|\cdot4\delta)2\pi \int\limits_{0}^{4\delta}\rho\log\frac{1}{\rho}d\rho =
o(1)\ as\ \delta\to0.
\end{align*}

\begin{align*}
I_6\ls&\iint\limits_{\Omega_2^+}|\log|\zeta-z|-\log|\zeta-z_0||
\sum\limits_{k=0}^{\infty} \frac{|\alpha|^{2k+1}|\zeta-z|^{2k}}{2^{2k}k!(k+1)!}
\times\\&\ \ \times (2(k+1) +
|\alpha|\cdot|\zeta-z|)d\xi\,d\eta+\\
 +& \iint\limits_{\Omega_2^+}|\log|\zeta-z_0||
\sum\limits_{k=0}^{\infty} \frac{|\alpha|^{2k+1}}{2^{2k}k!(k+1)!}
\left|(\zeta-z)^{2k}-(\zeta-z_0)^{2k}\right|\times\\
&\ \ \times \left|2(k+1) +
\alpha(\zeta-z)\right|d\xi\,d\eta+\\
+& \iint\limits_{\Omega_2^+}|\log|\zeta-z_0||
\sum\limits_{k=0}^{\infty} \frac{|\alpha|^{2k+1}}{2^{2k}k!(k+1)!}
|(\zeta-z_0)^{2k}|\cdot|\alpha|\cdot|z_0-z|d\xi\,d\eta.
\end{align*}

Using the inequalities
$$\left|\log\frac{|\zeta-z|}{|\zeta-z_0|}\right|<\frac{2\delta}{|\zeta-z_0|}
,\qquad \zeta\in \Omega_2^+,$$
$$\left|(\zeta-z)^{2k}-(\zeta-z_0)^{2k}\right|
\ls2k\left(C(\Gamma,z_0)\right)^{2k-1}\delta,\qquad \zeta\in
\Omega_2^+,$$ we have $$I_6\ls C(\Gamma,z_0,\alpha)\delta, $$
where $C(\cdot)$ denotes a constant depending only on the
parameters in the parenthesis.
\end{proof}

\stepcounter{subsection}
\begin{proof}[{\rm\thesubsection.} Proof of Theorem \ref{mait}]
Let us prove (\ref{phi+}) (the relation (\ref{phi-}) is
proved similarly).
We consider a sequence $z_n\in \Omega^+$, $z_n\to t\in \Gamma$, and denote by $%
\zeta_n$ the nearest to $z_n$ point of the curve $\Gamma$.

Applying formula (\ref{fi}) we have that
\begin{equation}  \label{fifi}
\begin{split}
&|\Phi_\alpha[f](z_n)-(I_{\alpha,\Gamma}(t)+1)f(t)-F_\alpha[f](t)| =\\
=&|\Phi_\alpha[f](z_n) - \Phi_\alpha[f(\zeta_n)](z_n) + \Phi_\alpha[%
f(\zeta_n)](z_n) - F_\alpha[f](\zeta_n) + F_\alpha[f](\zeta_n)- \\
&\hspace{52pt} - (I_{\alpha,\Gamma}(t)+1)f(t) - F_\alpha[f](t)|\ls M_1
+ M_2,
\end{split}
\end{equation}
where
\begin{align*}
M_1=&|\Phi_\alpha[f-f(\zeta_n)](z_n) - F_\alpha[f](\zeta_n)|, \\
M_2=&|(I_{\alpha,\Gamma}(z_n)+1)f(\zeta_n) + F_\alpha[f](\zeta_n) -
(I_{\alpha,\Gamma}(t)+1)f(t) - F_\alpha[f](t)|.
\end{align*}

Let $\alpha\not=0$. On the basis of relations (\ref{ha0}) --
(\ref{kal}) we have the representation
\begin{equation}\label{s}
    K_\alpha(z)=S_\alpha(z) + \varphi_\alpha(z),
\end{equation}
where $\varphi_\alpha(z)$ is a continuous function in $\C$ and
$$S_\alpha(z):=-\frac{1}{2\pi}\left(\frac{z}{|z|^2}+\alpha\log|z|\right).$$
Then
\begin{align*}
M_1\ls&\left|\int\limits_\Gamma S_\alpha(\zeta-z_n)\cdot\sigma\cdot (
    f(\zeta)-f(\zeta_n)) -\int\limits_\Gamma S_\alpha(\zeta-\zeta_n)
    \cdot\sigma\cdot(f(\zeta)-f(\zeta_n))\right|+\\
    +&\left|\int\limits_\Gamma \varphi_\alpha
    (\zeta-z_n)\cdot\sigma\cdot (
    f(\zeta)-f(\zeta_n)) -
    \int\limits_\Gamma \varphi_\alpha (\zeta-\zeta_n)
    \cdot\sigma\cdot(f(\zeta)-f(\zeta_n))\right|=:\\=:&M_3+M_4.
\end{align*}

By virtue of continuity of the functions $F_\alpha[f]$ (Lemma \ref{l2}), $%
I_{\alpha,\Gamma}$ (Lemma \ref{l4}), $\varphi_\alpha$ and $f$ we get that $M_2\to0$
and $M_4\to0$, when $z_n\to t$.

Let us fix an arbitrary $\varepsilon>0$. For any given $\delta>0$ we have

\begin{equation}
M_3\ls M_5+M_6+M_7,  \label{M3}
\end{equation}
where
\begin{align*}
M_5=&\left|\; \int\limits_{\Gamma_{\zeta_n,\delta}} S_\alpha(\zeta-\zeta_n)
\cdot\sigma\cdot (f(\zeta)-f(\zeta_n)) \right|, \\
M_6=&\left|\; \int\limits_{\Gamma_{\zeta_n,\delta}}S_\alpha(\zeta-z_n)
\cdot\sigma\cdot (f(\zeta)-f(\zeta_n)) \right|, \\
M_7=&\left|\; \int\limits_{\Gamma\setminus\Gamma_{\zeta_n,\delta}}
\left(S_\alpha(\zeta-z_n) - S_\alpha(\zeta-\zeta_n)\right) \cdot\sigma\cdot
(f(\zeta)-f(\zeta_n)) \right|.
\end{align*}

By  virtue  of  the equality (\ref{s}) it follows from the uniform
existence of $F_\alpha[f]$ (Lemma \ref{l2}) that for all
sufficiently small $\delta$ and for all $\zeta_n\in\Gamma$ the
inequality $M_5<\frac{\varepsilon}{3}$ is valid.

Let us estimate $M_6$. For any $\delta>0$ let us take $z_n$ near to $t$ so that
$|\zeta_n-z_n|<\frac\delta3$. We have
\begin{equation}
\begin{split}
M_6\ls&\left|\; \int\limits_{\Gamma_{\zeta_n,3|\zeta_n-z_n|}}
S_\alpha(\zeta-z_n) \cdot\sigma\cdot (f(\zeta)-f(\zeta_n)) \right|+ \\
+&\left|\; \int\limits_{\Gamma_{\zeta_n,\delta}\setminus
\Gamma_{\zeta_n,3|\zeta_n-z_n|}} S_\alpha(\zeta-z_n) \cdot\sigma\cdot
(f(\zeta)-f(\zeta_n)) \right|=: M_8+M_9.
\end{split}
\label{M4}
\end{equation}

Let us estimate $M_8$. Using the inequalities $|\zeta-\zeta_n|\ls3
|\zeta-z_n|<4\delta$, we obtain for sufficiently small $\delta<\frac{3}{4}$:

\begin{equation*}\label{ss}
\begin{split}
    \frac{|S_\alpha (\zeta-z_n)|}{|S_\alpha (\zeta -\zeta _n)|}=&
    \frac{\left|\frac{\zeta-z_n}{|\zeta-z_n|^2}+\alpha\log|\zeta-z_n|\right|}
    {\left|\frac{\zeta-\zeta_n}{|\zeta-\zeta_n|^2}+\alpha\log|\zeta-\zeta_n|\right|}\ls\\
\ls&\frac{\frac{3}{|\zeta-\zeta_n|}+|\alpha|\cdot|\log|\zeta-\zeta_n||+|\alpha|\log3}
    {\frac{1}{|\zeta-\zeta_n|}-|\alpha|\cdot|\log|\zeta-\zeta_n||}\ls 4,
\end{split}
\end{equation*}
\begin{equation}
|S_\alpha (\zeta -z_n)|=|S_\alpha (\zeta -\zeta _n)|\cdot \frac{|S_\alpha (\zeta
-z_n)|}{|S_\alpha (\zeta -\zeta _n)|}\ls4\cdot|S_\alpha (\zeta -\zeta _n)|.
\label{K1}
\end{equation}

Due to the uniform existence of the integral (\ref{psi}) and the equality (\ref{s})
it follows from (\ref{K1}) that for all sufficiently small $\delta$ and for all
$|\zeta_n-z_n|<\frac{\delta}{3}$
\begin{equation}
M_8\ls 4\int\limits_{\Gamma_{\zeta_n,3|\zeta_n-z_n|}} |S_\alpha \left( \zeta -\zeta
_n\right) |\cdot |\sigma |\cdot |f(\zeta )-f( \zeta _n) |<\frac \varepsilon 6.
\label{M6}
\end{equation}

Let us estimate $M_9$. We get
\begin{equation}
|S_\alpha(\zeta-z_n)|\ls|S_\alpha(\zeta-\zeta_n)| +
|S_\alpha(\zeta-z_n)-S_\alpha(\zeta-\zeta_n)|.  \label{M71}
\end{equation}

As long as $|z_n-\zeta_n|\ls\frac12|\zeta-z_n|$,
$3|\zeta_n-z_n|<|\zeta-\zeta_n|\ls\delta<\frac{3}{4}$ and
$\frac{1}{|\zeta-\zeta_n|}\ls2\pi|S_\alpha(\zeta-\zeta_n)|$ we have
\begin{equation}
\begin{split}
|S_\alpha(\zeta-z_n)-&S_\alpha(\zeta-\zeta_n)| \ls \\
\ls&\frac{1}{2\pi} \left|\frac{\zeta-z_n}{
 |\zeta-z_n|^2}- \frac{\zeta-\zeta_n}{|\zeta-\zeta_n|^2}\right| +
 \frac{|\alpha|}{2\pi}\left|\log\frac{|\zeta-z_n|}{
 |\zeta-\zeta_n|}\right|=\\
 =&\frac{|z_n-\zeta_n|}{2\pi\,
 |\zeta-z_n|\cdot|\zeta-\zeta_n|} +
 \frac{|\alpha|}{2\pi}\left|\log\frac{|\zeta-z_n|}{
 |\zeta-\zeta_n|}\right|\ls\\
 \ls&\frac{1}{4\pi\, |\zeta-\zeta_n|} +
 \frac{|\alpha|}{2\pi}\log\frac{3}{2}\ls \frac{1+|\alpha|}{4\pi\, |\zeta-\zeta_n|}\ls
 \frac{1+|\alpha|}{2}|S_\alpha(\zeta-\zeta_n)|.
  \label{M72}
\end{split}
\end{equation}

From (\ref{M71}), (\ref{M72}) we get
\begin{equation}
\begin{split}
M_9\ls&\frac{3+|\alpha|}{2}\int\limits_{\Gamma_{\zeta_n,\delta}\setminus
\Gamma_{\zeta_n,3|\zeta_n-z_n|}} |S_\alpha(\zeta-z_n)| \cdot|\sigma|\cdot
|f(\zeta)-f(\zeta_n)|<\frac{\varepsilon}{6}
\end{split}
\label{M78}
\end{equation}
for sufficiently small $\delta$ and for $|\zeta_n-z_n|<\frac{\delta}{3}$.

So we have from (\ref{M4}), (\ref{M6}, (\ref{M78}) that $M_6<\frac{\varepsilon}{3}$.

In order to estimate $M_7$ fix any $\delta$ satisfying  all the
conditions stated above and take $z_n$ near  $t$ so that
$|\zeta_n-z_n|\ls\frac\delta3$.

We have $\delta<|\zeta-\zeta_n|$, $\frac23\delta<|\zeta-z_n|$. Therefore
\begin{equation}
\frac{|z_n-\zeta_n|}{|\zeta-z_n|\cdot|\zeta- \zeta_n|}\ls
\frac{3}{2\delta^2}|z_n-\zeta_n|,  \label{M60}
\end{equation}
and by  Lagrange's theorem $$\left|\log\frac{|\zeta-z_n|}{|\zeta
-\zeta_n|}\right|=\frac1\mu\,||\zeta-z_n| -|\zeta -\zeta_n||\ls
\frac3{2\delta} |z_n-\zeta_n|,$$ where $\mu$ lies between
$|\zeta-z_n|$ and $|\zeta -\zeta_n|$. That is why using the
relation (\ref{M72}), we get $$
|S_\alpha(\zeta-z_n)-S_\alpha(\zeta-\zeta_n)|\ls
3\frac{1+|\alpha|\delta}{4\pi\delta^2}|z_n-\zeta_n| ,$$ and taking
into account the boundedness of the function $f$, we obtain for
the above fixed $\delta$ and for $z_n$ sufficiently near to $t$
\begin{equation*}
\begin{split}
M_7\ls& \int\limits_{\Gamma\setminus\Gamma_{\zeta_n,\delta}} |S_\alpha(\zeta-z_n)
-S_\alpha(\zeta-\zeta_n)| \cdot|\sigma|\cdot |f(\zeta)-f(\zeta_n)|\ls\\
\ls&\frac{1+|\alpha|\delta}{2\delta^2}\,l(\Gamma)\,\max\limits_{t\in\Gamma}|f(t)|\,
|\zeta_n-z_n| <\frac\ve3,
\end{split}
\end{equation*}
where $l(\Gamma)$ denotes the length of $\Gamma$.

Thus, we have $M_3<\ve$ and, consequently, the relation (\ref{phi+}) is proved.

The continuity of $\Phi^\pm_\alpha[f]$ on $\Gamma$ follows now
from Lemmas \ref{l2} and \ref{l4}. This completes the proof of
Theorem \ref{mait}.
\end{proof}

\section{Proof of main results}

\subsection{}We identify a complex quaternion $a=\sum_{k=0}^3a_k\bs i_k$ with the
scalar-vector pair $(a_0,\bs a)$, where $\bs a=\sum_{k=1}^3a_k\bs
i_k$ is a vector of the complex linear space $\C^3$ with the
canonical basis $\bs i_1$, $\bs i_2$, $\bs i_3$. Then a
quaternionic function $f=\sum_{k=0}^3f_k{\bs i_k}$ is
interpretable as a pair $\mathcal F=(f_0,\bs f)$, operator
$\partial_\alpha$  as a pair $(M^\alpha,\,_{st}\bs\partial)$,
where $_{st}\bs\partial:=\partial_1\bs i_1+\partial_2\bs i_2$.
Using the vectorial  representation of the multiplication of any
complex quaternions $a=(a_0,\bs a)$ and $b=(b_0,\bs b)$ (see
\cite{krs}, p. 24):
\begin{equation}\label{qv}
    ab=(a_0b_0-<\bs a,\bs b>,[\bs a,\bs b]+a_0\bs b+b_0\bs a),
\end{equation}
we obtain
\begin{equation*}
\partial_{\alpha}f =(\alpha f_0-div\bs f,\bs{rot f}+ \alpha\bs f +
\bs{grad}\,f_0),
\end{equation*}
reducing to the system (\ref{system2}) as the vector form of
Definition \ref{df} of an $\alpha$-hy\-per\-ho\-lo\-mor\-phic
function.

\stepcounter{subsection}
\begin{proof}[{\rm\thesubsection.} Proof of Theorem \ref{theo6}] The representation
(\ref{kalv}) follows from the formula (\ref{kal}) and we obtain
(\ref{fiav}) from (\ref{fia}) by using the equality (\ref{qv}).
Combining the vector form of the functions $F_\alpha$,
$I_{\alpha,\Gamma}$, $\Phi_\alpha^\pm$ in Theorem \ref{mait} with
the equality (\ref{qv}), we arrive  at Theorem \ref{theo6} as a
vector reformulation of Theorem \ref{mait}.
\end{proof}
\stepcounter{subsection}\begin{proof}[{\rm\thesubsection.} Proof of Corollary
\ref{17}]
Applying Theorem \ref{theo6} to the pair $\mathcal{F}=(0,\bs f)$, we
obtain the desired conclusion. Because of the condition $\bs f\in\mathcal
M(\Gamma;\C^3)$ the Cauchy-type integral $\Phi_\alpha[\mathcal F]$ is purely
vectorial and therefore its boundary values $\Phi_\alpha^\pm[\mathcal F]$ are also
purely vectorial.
\end{proof}

\section*{Acknowledgments}
The first-named author was supported in part by INTAS-99-00089 and by \linebreak[4]
CONACYT project.

The second-named author was partially supported by CONACYT projects as well as by
Insituto Polit\'ecnico Nacional in the framework of COFAA and CGPI programs.

\end{document}